\documentclass[amssymb,amstex,amsmath,11pt]{article}
\usepackage[margin=1.4in]{geometry}
\usepackage{amsmath,amssymb,latexsym,amsfonts,url}
\usepackage{cite}
\usepackage{graphicx}
\usepackage{amsthm}
\usepackage{tikz}
\usepackage{pgfplots}
\usepackage{color}
\usepackage{epstopdf}
\usetikzlibrary{intersections, pgfplots.fillbetween}
\usepackage[all]{xy}
\pgfplotsset{compat=newest}
\pgfdeclarelayer{pre main}

\numberwithin{equation}{section}

\renewcommand{\epsilon}{\varepsilon}
\theoremstyle{plain}
\newtheorem{thm}{Theorem}[section]

\newtheorem{remark}[thm]{Remark}
\newtheorem{p}[thm]{Problem}

\newtheorem*{theorem*}{Theorem}
\newtheorem*{proposition*}{Proposition}
\theoremstyle{definition}
\newtheorem{df}[thm]{Definition}

\theoremstyle{remark}

\newcommand{\dist}{{\rm dist}}

\newcommand{\tr}{{\rm{tr}}}

\newcommand{\C}{\mathbb{C}}

\def\SS{{\mathbb S}}

\def\RR{{\mathbb R}}
\def\C{{\mathcal{C}}}

\def\({\left(}
\def\){\right)}
\def\[{\left[}
\def\]{\right]}
\def\<{\left\langle}
\def\>{\right\rangle}

\def\pa {\partial}
\def\hess {\mathrm{Hess}\,}
\def\tr {\mathrm{tr}}
\def\ric{\mathrm{Ric}}
\def\Om{\Omega}

\def\grad{\nabla}
\def\trhess{\tr(\hess^2u)}

\def\ga{\Gamma_0 } 
\def\gacone{\Gamma_1} 
\def\p{\widetilde{P}}
\def\f{h'}

\def\RR{{\mathbb R}}

\def\Om{\Omega}

\def\pa{\partial }

\def\({\left(}
\def\){\right)}

\begin{document}
\begin{title}
{Radial symmetry and partially overdetermined problems in a convex cone}
\end{title}
\begin{author}{Jihye Lee \and Keomkyo Seo}\end{author}


\maketitle

\begin{abstract}
\noindent We obtain the radial symmetry of the solution to a partially overdetermined boundary value problem in a convex cone in space forms by using the maximum principle for a suitable subharmonic function $P$ and integral identities. In dimension $2$, we prove Serrin-type results for partially overdetermined problems outside a convex cone. Furthermore, we obtain a Rellich identity for an eigenvalue problem with mixed boundary conditions in a cone. \\


\end{abstract}

\section{Introduction}
In a celebrated paper \cite{Serrin}, James Serrin obtained the following remarkable result. Let $\Omega$ be a smooth, bounded, open, connected domain in $\mathbb{R}^n$ and let $\nu$ be the outward unit normal to $\partial \Omega$. If $u$ is a smooth solution to the overdetermined problem
\begin{eqnarray*}
\begin{cases}
\Delta u =-1 & \mbox{in } \Omega, \\
u=0 & \mbox{on } \partial  \Omega, \\
\frac{\partial u}{\partial \nu} = const=c & \mbox{on } \partial \Omega,
\end{cases}
\end{eqnarray*}
then $\Omega$ is a ball and the solution $u$ is radially symmetric. His proof is based on the moving plane method (or Alexandrov \cite{Alexandrov} reflection method). After Serrin's work, his symmetry result has been generalized to space forms (see \cite{CV2017, CV2019, FMW, KP, Molzon, QX, Roncoroni} for instance and references therein).

On the other hand, Pacella-Tralli \cite{PT} studied a partially overdetermined problem for a domain in a convex cone in the Euclidean space. In order to describe precisely, let us introduce some notations. Let $\mathcal{C}$ be an open convex cone with vertex at the origin $O$ in $\mathbb{R}^n$, $n\geq 2$, i.e.,
 $$\C = \{tx:x\in \omega, t\in (0,\infty) \}$$
 for some domain $\omega$ in the unit sphere $\mathbb{S}^{n-1}$. Recall that a $C^2$ domain $\Omega$ is {\it convex} in an $n$-dimensional Riemannian manifold $M^n$ if the normal curvature on $\partial \Omega$ with respect to the inward normal direction is nonnegative. In particular, if the domain is convex and $n\geq3$, the second fundamental form II is nonnegative at every point on the boundary $\partial \Omega$. Given an open convex cone $\mathcal{C}$ such that $\partial \mathcal{C} \setminus \{O\}$ is smooth and a domain $\Omega \subset \mathcal{C}$, we denote by $\Gamma_0$ its relative boundary, i.e., $\Gamma_0 = \partial \Omega \cap \mathcal{C}$ and let $\Gamma_1=\partial \Omega \setminus \overline{\Gamma_0}$. Assume that $\mathcal{H}^{n-1} (\Gamma_i)>0$ for $i=0,1$, where $\mathcal{H}^{n-1}$ denotes the $(n-1)$-dimensional Hausdorff measure. Moreover, we assume that $\Gamma_0$ is an $(n-1)$-dimensional smooth manifold, while $\partial \Gamma_0 = \partial \Gamma_1 \subset \partial \mathcal{C} \setminus \{O\}$ is an $(n-2)$-dimensional smooth manifold. Following \cite{PT}, such a domain $\Omega$ is called a {\it{sector-like}} domain. Given a sector-like domain $\Omega$ in an open convex cone $\mathcal{C}\subset \mathbb{R}^n$, consider the partially overdetermined mixed boundary value problem
\begin{eqnarray}
\begin{cases}
\Delta u = -1 & \mbox{in } \Om, \\
u =0, \quad \frac{\pa u}{\pa \nu} = const=-c<0 & \mbox{on } \ga,\\  \label{problem: partially overdetermined problem}
        \frac{\pa u}{\pa \nu} = 0 & \mbox{on } \gacone \setminus \{O\},
\end{cases}
\end{eqnarray}
where $\nu=\nu(x)$ denotes the outward unit normal to $\partial \Omega$ wherever is defined (that is for $x\in \Gamma_0\cup \Gamma_1 \setminus \{O\}$). Pacella-Tralli \cite{PT} proved the following.
\begin{theorem*}[\cite{PT}]
Let $\Omega$ be a sector-like domain in an open convex cone $\mathcal{C}$ in $\mathbb{R}^n$. Assume that there exists a classical solution $u \in C^2 (\Omega)\cap C^1 (\Omega \cup \Gamma_0 \cup \Gamma_1 \setminus \{O\})$ to the partially overdetermined problem (\ref{problem: partially overdetermined problem}) such that $u\in W^{1,\infty} (\Omega) \cap W^{2,2} (\Omega)$. Then
$$\Omega=\mathcal{C} \cap B_R (p) \quad \mbox{and} \quad  u(x)=\frac{R^2-|x-p|^2}{2n}$$
for some point $p \in \partial \mathcal{C}$. Here $B_R(p)$ denotes the ball centered at a point $p\in \mathbb{R}^n$ and $R=nc$.
\end{theorem*}
\noindent Note that the point $p$ may not be the origin $O$ in the above theorem. Recently, Ciraolo-Roncoroni \cite{CR} extended the above theorem into space forms. Indeed, they considered the partially overdetermined problem in space forms
\begin{eqnarray}
\begin{cases}
\Delta u +nKu= -n & \mbox{in } \Om, \\
u =0, \quad \frac{\pa u}{\pa \nu} = const=-c & \mbox{on } \ga,\\  \label{problem S: partially overdetermined problem in S}
        \frac{\pa u}{\pa \nu} = 0 & \mbox{on } \gacone \setminus \{O\},
\end{cases}
\end{eqnarray}
where $K=0, +1, -1$ in the Euclidean space, in the upper unit hemisphere $\mathbb{S}^n_+$, and in the hyperbolic space $\mathbb{H}^n$, respectively. In \cite{CR}, Ciraolo-Roncoroni obtained  that if $\Omega$ is a sector-like domain in a convex cone in space forms, then $\Omega=\mathcal{C} \cap B_R (p)$ and the solution $u$ is radially symmetric with respect to the point $p$, where $B_R (p)$ denotes a geodesic ball of radius $R$ centered at $p$.

This paper is organized as follows. In Section $2$, we investigate two partially overdetermined problems for domains on a convex cone with vertex at $p$ in the unit sphere $\mathbb{S}^n$. Firstly, we consider the equation
$$\Delta u = - n \cos r,$$
where $\Delta$ denotes the Laplace-Beltrami operator on $\mathbb{S}^n$ and $r$ denotes the geodesic distance from the vertex $p$. In fact, using the above equation, Molzon \cite{Molzon} extended Serrin's symmetry result to the upper unit hemisphere $\mathbb{S}^n_+$. In Theorem \ref{thm1}, we obtain an analogue of Molzon's result for domains in a convex cone. Secondly, we consider the same partially overdetermined problem (\ref{problem S: partially overdetermined problem in S}) for a domain in a cone in the unit sphere as in \cite{CR}. However, we do not assume that the domain $\Omega$ is contained in the upper hemisphere $\mathbb{S}^{n}_+$, but assume that $\Omega$ is star-shaped with respect to the vertex $p$. Using the maximum principle for a suitable subharmonic function $P$ and some integral identities (which is originated by Weinberger \cite{Weinberger}), we are able to prove a rigidity result of Serrin type for a star-shaped domain in a convex cone with vertex at $p$ in Theorem \ref{thm2}.

In $2$-dimensional case, it turns out that it is not necessary to assume that a domain is contained in a convex cone for the partially overdetermined problem (\ref{problem S: partially overdetermined problem in S}). Indeed, we prove that if $\mathcal{C} \subset \mathbb{S}^2$ is a convex cone with vertex at $p$ and $\Omega$ is a star-shaped domain with respect to $p$ {\it outside} a convex cone $\mathcal{C}$ and if (\ref{problem S: partially overdetermined problem in S}) admits a solution, then $\Omega = B_R (p_0) \setminus \overline{\mathcal{C}}$ for some $p_0\in \partial \mathcal{C}$ and the solution $u$ is radially symmetric in Section $3$ (see Theorem \ref{thm3.2}). In Section $4$, we study an eigenvalue problem with mixed boundary conditions in a cone. Given the Dirichlet eigenvalue problem for a bounded domain $\Omega \subset \mathbb{R}^n$

\begin{eqnarray*}
\begin{cases}
\Delta u +\lambda u = 0 & \mbox{in }  \Om, \\
u =0  & \mbox{on } \partial \Omega,
\end{cases}
\end{eqnarray*}
\noindent Rellich \cite{Rellich} obtained the following identity
\begin{eqnarray}
\lambda = \frac{1}{4} \int_{\partial \Omega} \left(\frac{\partial u}{\partial \nu}\right)^2 \frac{\partial r^2}{\partial \nu} d\sigma, \label{identity: Rellich}
\end{eqnarray}
where $r$ denotes the distance from the origin and $u$ is normalized so that $\int_\Omega u^2 dV =1$. (\ref{identity: Rellich}) is called the {\it Rellich identity}. In 1991, Molzon \cite{Molzon} extended (\ref{identity: Rellich}) to space forms. Motivated by his result, we consider the mixed boundary eigenvalue problem for a domain $\Omega$ in a cone $\mathcal{C}$ with vertex at $p$

\begin{eqnarray*}
\begin{cases}
            \Delta u +\lambda u = 0 & \mbox{in } \Omega, \\
            u=0 & \mbox{on } \partial \Omega \cap \mathcal{C}, \\
            \frac{\pa u}{\pa \nu} = 0 & \mbox{on } \partial \Omega \setminus (\{p\} \cup (\overline{\partial \Omega \cap \mathcal{C}})).
            \end{cases}
\end{eqnarray*}
\noindent and obtain a similar result for domains in a cone  in Theorem \ref{thm4.1}. Note that the cone $\mathcal{C}$  needs not to be convex.

\section{Partially overdetermined problems inside a convex cone in the unit sphere}
Let $(M^n, g)$ be an $n$-dimensional space form, i.e., an $n$-dimensional complete simply-connected Riemannian manifold with constant sectional curvature $K$. Up to homotheties, we may assume that $K=0, 1$, and $-1$: the corresponding spaces are the Euclidean space $\mathbb{R}^n$, the unit upper hemisphere $\mathbb{S}^n_+$, and the hyperbolic space $\mathbb{H}^n$, respectively. These three spaces can be represented as the warped product space $M= I \times \mathbb{S}^{n-1}$ which is equipped with the rotationally symmetric metric
$$g=dr^2 + h(r)^2 g_{\mathbb{S}^{n-1}},$$
where $g_{\mathbb{S}^{n-1}}$ denotes the round metric on the $(n-1)$-dimensional unit sphere $\mathbb{S}^{n-1}$ and
\begin{itemize}
\item $h(r) = r$ on $I = [0,\infty)$ in $\mathbb{R}^n$;
\item $h(r) = \sin r$ on $I = [0,\frac{\pi}{2})$ in $\mathbb{S}^n_+$;
\item $h(r) = \sinh r$ on $I = [0,\infty)$ in $\mathbb{H}^n$.
\end{itemize}
Here $r(\cdot)$ denotes the distance $\dist (\cdot, p)$ from the {\it pole} $p$ of the model space. Define a {\it cone} $\mathcal{C}$ with vertex at $p$ as follows:
$$\mathcal{C}:= \{tx: x\in \omega, t \in I\}$$
for some domain $\omega \subset \mathbb{S}^{n-1}$. Note that $\mathcal{C} \subset M$ is convex if $\omega$ is convex in $\mathbb{S}^{n-1}$.
\begin{df}
A connected bounded open set $\Om \subset \C$ is an \textit{admissible interior} domain if the boundary $\pa \Om$ satisfies the following.
\begin{itemize}
\item $\pa \Om$ contains the vertex $p$.
\item $\ga:=\pa \Om \setminus \pa \C \neq \emptyset$ is an $(n-1)$-dimensional smooth manifold.
\item $\gacone:= \pa \Om \setminus \overline{\ga} \neq \emptyset$ and $\partial \Gamma_0 = \partial \Gamma_1 \subset \partial \mathcal{C} \setminus \{p\}$ is an $(n-2)$-dimensional smooth manifold.
\item $\mathcal{H}^{n-1} (\Gamma_i) >0$ for $i=0,1$, where $\mathcal{H}^{n-1}$ denotes the $(n-1)$-dimensional Hausdorff measure.
\end{itemize}
\end{df}

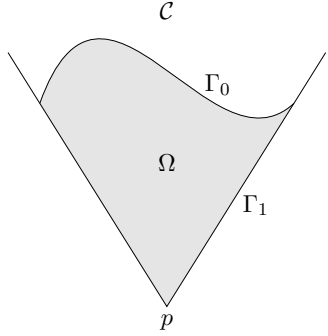
\begin{figure}[h!t!b]
\begin{center}
\resizebox{0.3\textwidth}{!}{%
\begin{tikzpicture}
\pgfsetlayers{pre main,main}
\draw[name path=A] (2.5,4)--(0,0)--(-2.5,4);
\draw[name path=C] (-2, 3.2) .. controls (-1, 6) and (0.7, 2) .. (2, 3.2);
\tikzfillbetween[of=C and A]{gray, opacity=0.2};
\coordinate[label={$p$}] (p) at (0,-0.5);
\coordinate[label={$\Om$}] (Om) at (0,2);
\coordinate[label={$\gacone$}] (gacone) at (1.4,1.3);
\coordinate[label={$\ga$}] (ga) at (0.8,3.2);
\coordinate[label={$\C$}] (C) at (0,4.4);
\end{tikzpicture}
}
\caption{An admissible interior domain $\Omega$ inside a convex cone $\mathcal{C}$}\label{figure: interior}
\end{center}
\end{figure}
\noindent Following \cite{CR, PT}, we remark that if the boundary of a sector-like domain contains the vertex $p$, then such a domain is an admissible interior domain. Modifying Molzon's argument in \cite{Molzon}, we are able to prove the following partially overdetermined problem for domains in a convex cone in the upper unit hemisphere $\mathbb{S}^n_+$.
\begin{thm}\label{thm1}
Let $\C \subset \SS_+^n$ be an open convex cone with  vertex at $p$ such that $\pa \C \setminus \{p\}$ is smooth and $\Om$ be an admissible interior domain in $\mathcal{C}$.  Suppose there exists a solution $u\in C^2(\overline{\Om})$ satisfying
            $$\begin{cases}
                \Delta u = -n h'(r)= -n \cos r  & \mbox{in } \Om, \\
                u =0, \quad \frac{\pa u}{\pa \nu} = const=c & \mbox{on } \ga,\\
                \frac{\pa u}{\pa \nu} = 0 & \mbox{on } \gacone \setminus\{p\},
            \end{cases}$$
where $\nu$ is the outward unit normal to $\ga \cup \gacone \setminus \{p\}$ and  $r(x) = \mathrm{dist}(p,x)$. Then $\Om$ is part of the ball centered at the vertex $p$ of the cone $\mathcal{C}$, i.e.,
            $$\Om = \C\cap B_R (p),$$
        where $B_R (p)$ denotes the geodesic ball centered at the vertex $p$  with radius $R$ in the upper hemisphere $\SS_+^n$. 
        Moreover, the solution $u$ is radial and it is given by
            $$u(x) = \cos r - \cos R,$$
        where $R = \sin^{-1}(-c)$.
\end{thm}

    \begin{proof}
       We first claim that $u>0$ in $\Om$. To see this, let
       $$u^{-}(x) :=\begin{cases}
                0 & \mbox{if } \, u(x) \geq 0, \\
                u(x) & \mbox{if } \, u(x) < 0.
       \end{cases}$$
        Then $u^{-} (x) \leq 0$ and $\Delta u = - n \cos r \leq 0$ on $\Omega \subset \SS_+^n$.
        Using the divergence theorem, we obtain
        \begin{align*}
          0 & \leq \int_{\Om} u^{-} \Delta u dV \\
           & = \int_{\pa \Om} u^{-} \frac{\pa u}{\pa \nu} d \sigma - \int_{\Om} \<\grad u^{-} , \grad u\> dV\\
           & = - \int_{\Om \cap \{ u<0\}} |\grad u|^2 dV \leq 0,
        \end{align*}
        which yields $u \geq 0$ in $\Omega$. Since $u$ is not constant, it follows from the maximum principle that $u>0$ in $\Omega$. Thus we may assume that $u$ is positive in $\Omega$.

        It is well-known that
            \begin{equation*}
                (\Delta u)^2 \leq  n \,\trhess,
            \end{equation*}
        where $\hess u$  denotes the Hessian of $u$ and $\hess^2 u = \hess u \circ \hess u$. Thus
            \begin{equation}\label{ineq::schwarz}
                nh'^2 \leq \trhess.
            \end{equation}
        A straightforward calculation shows that
        \begin{equation*}
          \hess \f = - \f g \quad \mbox{ and } \quad \Delta \f = -n\f,
        \end{equation*}
        where $g$ is the standard metric of $\SS^n$.  Applying the polarized Bochner formula on a Riemannian manifold
\begin{equation*}
\Delta \langle \nabla  \phi, \nabla \psi \rangle = \nabla \phi (\Delta \psi) + \nabla \psi (\Delta \phi) + 2 \tr (\hess \phi \circ \hess \psi) + 2\ric (\nabla \phi, \nabla \psi)
\end{equation*}
for any smooth function $\phi$ and $\psi$, we have
            \begin{align}\label{eq::boch}
                \Delta \< \grad (u-\f), \grad u \>   &= \< \grad (\Delta (u-\f) ), \grad u \> + \< \grad (u-\f), \grad (\Delta u)\> \nonumber               \\
                & \quad + 2 \tr (\hess (u-\f) \circ \hess u) + 2 \ric (\grad (u-\f), \grad u),
            \end{align}
        where $\ric(\cdot,\cdot)$ is the Ricci tensor of $g$.
        Using (\ref{ineq::schwarz}), we get
            $$\mathrm{tr}(\hess (u-\f) \circ \hess u) = \trhess +\f \Delta u  = \trhess -nh'^2 \geq 0.$$
        Thus (\ref{eq::boch}) becomes
            $$\Delta \< \grad (u-\f ), \grad u\> \geq  -n  \< \grad (u-\f), \grad \f\> + 2(n-1) \< \grad (u-\f), \grad u\>.$$
        Since $u>0$ in $\Om$, we obtain
            \begin{align}\label{ineq::boch}
                \int_{\Om} u \Delta \<\grad (u-\f), \grad u\> dV    &\geq -n \int_{\Om} u \<\grad (u-\f), \grad \f \> dV \nonumber      \\
                &\quad + 2 (n-1) \int_{\Om} u \<\grad (u-\f), \grad u\> dV.
            \end{align}
        Note that
        $$\frac{\partial h'}{\partial \nu} = \< \nabla h', \nu \> = -\sin r \< \grad r, \nu\> = 0$$
        on $\Gamma_1$.
        Using the divergence theorem, we get
            \begin{align*}
                \nonumber\int_{\Om} \< \grad (u-\f), \grad (u^2) \> dV & = \int_{\Om} \mathrm{div}( u^2 \grad (u-\f) )dV - \int_{\Om} u^2 \Delta(u-\f) dV \\
                &= \int_{\pa \Om} u^2 \frac{\pa}{\pa \nu} (u-\f) d\sigma \\
                &= 0,
            \end{align*}
        which yields
            \begin{equation}\label{u0}
                \int_{\Om} u \<\grad (u-\f) ,\grad u\> dV = 0.
            \end{equation}
        From (\ref{u0}), we see that
            \begin{align}\label{eq::umf2}
                \nonumber\int_{\Om} u \< \grad (u-\f), \grad \f \> dV
                &= -\int_{\Om} u\<\grad (u-\f), -\grad \f\> dV \\
                &= - \int_{\Om} u |\grad (u-\f) |^2 dV.
            \end{align}
        Combining (\ref{ineq::boch}), (\ref{u0}), and (\ref{eq::umf2}), we get
            \begin{align}\label{geq}
            \int_{\Om} u \Delta \< \grad (u-\f), \grad u\> dV & \geq n \int_{\Om} u |\grad (u-\f) |^2dV\geq 0.
            \end{align}
        On the other hand, Green's identity gives
            \begin{align}\label{uf0}
              \int_{\Om} u \Delta \< \grad (u-\f), \grad u\> dV =& \int_{\Om} \< \grad (u-\f), \grad u\> \Delta u dV \nonumber\\
                & + \int_{\pa \Om} u \frac{\pa}{\pa \nu} \< \grad (u-\f), \grad u\> d\sigma \nonumber\\
                 & - \int_{\pa \Om} \< \grad (u-\f), \grad u \> \frac{\pa u}{\pa \nu} d \sigma.
            \end{align}
        Using the divergence theorem and the boundary conditions, we get
            \begin{align}
                \int_{\Om} \< \grad (u-\f) , \grad (u\f) \>dV & =  \int_{\Om} \mathrm{div}(u\f \grad(u-\f)) dV - \int_{\Om} u\f \Delta(u-\f) dV \nonumber\\
                & =  \int_{\pa \Om} u\f \frac{\pa }{\pa \nu} (u-\f) d\sigma  \nonumber\\
                & = 0. \label{ineq: 111}
            \end{align}
       Using \eqref{eq::umf2}  and \eqref{ineq: 111}, we have
            \begin{align}\label{uf}
              \int_{\Om} \<\grad (u-\f), \grad u \> \Delta u dV &= -n\int_{\Om}\f \< \grad (u-\f) , \grad u\> dV \nonumber\\
                & = n \int_{\Om} u \< \grad (u-\f) , \grad \f \> dV \nonumber\\
                & = -n \int_{\Om} u |\grad (u-\f)|^2 dV \nonumber\\
                & \leq 0.
            \end{align}
        Since $\Gamma_0 \subset \{u=0\}$ is a level set of $u$,
                $$ \grad u = c\nu \quad \mbox{ on }\Gamma_0.$$
        Thus
        \begin{align}\label{uf1}
          \int_{\pa \Om} \< \grad (u-\f), \grad u \> \frac{\pa u}{\pa \nu} d \sigma
           & = c^2\int_{\Gamma_0} \<\grad (u-\f), \nu\> d \sigma \nonumber\\
           &= c^2 \int_{\pa \Om} \< \grad (u-\f), \nu \> d \sigma \nonumber\\
           &= c^2 \int_{\Om} \Delta (u-\f) dV \nonumber \\
           &= 0.
        \end{align}
        Substituting (\ref{uf}) and (\ref{uf1}) into (\ref{uf0}), we obtain
            \begin{align}\label{temp1}
                \int_{\Om} u \Delta \< \nabla (u-\f), \grad u\> dV \leq
                \int_{\gacone}   u\frac{\pa}{\pa \nu} \< \grad (u-\f), \grad u\> d\sigma.
            \end{align}
        One can evaluate the right-hand side of (\ref{temp1}). We note that
            \begin{align}\label{gacone:in}
                \frac{\pa }{\pa \nu} \< \grad (u-\f),\grad u \> \nonumber
                & = \< \grad_\nu \grad(u-\f) , \grad u \> + \< \grad (u-\f), \grad_\nu \grad u\>  \\
                & = 2 \hess u (\grad u,\nu) - \hess \f ( \grad u, \nu ) - \hess u ( \grad \f, \nu)\nonumber\\
                & = 2 \hess u (\grad u,\nu) + \f \frac{\pa u}{\pa \nu} - \hess u ( \grad \f, \nu).
            \end{align}
        Since $\frac{\pa u}{\pa \nu} = 0$ and $\frac{\pa \f}{\pa \nu} = 0$ on $\gacone$, $\grad u$ and $\grad \f $ are tangent vectors on $\gacone$. Observe that
            $$\grad_{\grad \f} \nu = 0 \quad \mbox{ on } \gacone.$$
        On $\Gamma_1$,
            \begin{equation}\label{gradf}
                0 = \grad_{\grad \f} \(\frac{\pa u}{\pa \nu}\) = \< \grad_{\grad \f} \grad u, \nu \> + \< \grad u, \grad_{\grad \f} \nu \> = \hess u ( \grad h', \nu).
            \end{equation}
        Moreover, the convexity of the cone $C$ tells us that
        $$\mathrm{II}(\grad u,\grad u) = \< \nabla_{\nabla u} \nu, \nabla u\> \geq 0 \quad \mbox{ on } \gacone, $$
        where $\mathrm{II}(\cdot, \cdot)$ is the second fundamental form. Thus
            \begin{align}\label{II}
                0 = \grad_{\grad u} \(\frac{\pa u}{\pa \nu}\) &= \< \grad_{\grad u} \grad u, \nu \> + \< \grad u, \grad_{\grad u} \nu \> \nonumber\\
                &=\hess u (\grad u, \nu)+ \mathrm{II}(\grad u,\grad u) \nonumber \\
                &\geq \hess u (\grad u, \nu)  \quad \mbox{ on } \Gamma_1.
            \end{align}
        Plugging (\ref{gradf}) and   (\ref{II}) into  (\ref{gacone:in}), we get
            \begin{equation*}
                \frac{\pa }{\pa \nu} \< \grad (u-\f),\grad u \>  \leq 0 \quad \mbox{ on } \gacone.
            \end{equation*}
       By the continuity of $u$, we see that $u \geq 0$ on $\gacone$, since $u>0$ in $\Omega$. Thus
            \begin{equation}\label{int::gacone}
                \int_{\gacone} u \frac{\pa }{\pa \nu} \< \grad (u-\f), \grad u\> d\sigma \leq 0.
            \end{equation}
        Therefore (\ref{temp1}) and (\ref{int::gacone}) shows
            \begin{equation}\label{leq}
              \int_{\Om} u \Delta \< \grad (u-\f), \grad u\> dV\leq 0.
            \end{equation}
        By (\ref{geq}) and (\ref{leq}),
            $$\int_{\Om} u \Delta \< \grad (u-\f), \grad u\> dV = 0.$$
       Equality in (\ref{geq}) shows
        $$\grad (u-\f) \equiv 0 \mbox{ in } \Omega,$$ which implies that
            $$u(x) = \f + a = \cos r + a$$
        for some constant $a$.
       Moreover, the constant $a$ can be expressed  in terms of the constant $c$. To see this, note that the function $u$ vanishes on $\ga$  by the boundary condition. Thus $\ga$ is part of the boundary of the geodesic ball $B_R (p)$ of radius $R= \cos^{-1}(-a)$ centered at $p$. This shows that
            $$u(x) = \cos r - \cos R$$
        and
         $$\Omega = B_R (p) \cap \C.$$
        The boundary condition on $\ga$ gives
            $$c = \frac{\pa u}{\pa \nu} = \< -\sin R \,\grad r, \grad r\> = -\sin R \quad \mbox{ on } \ga.$$
       Finally we obtain
            $$u(x) = \cos r - \cos R,$$
        where $R= \sin^{-1}(-c).$
    \end{proof}

\begin{remark} {\rm
        Using the divergence theorem, we get
            $$-n \int_{\Om} h' dV = \int_{\Om} \Delta u \,dV = \int_{\pa \Om} \frac{\pa u}{\pa \nu}\, d\sigma = \int_{\ga} \frac{\pa u}{\pa \nu} \,d \sigma = c |\ga|.$$
        Since $h' \geq 0$ on $\Om \subset \SS_{+}^n$,
            $$c=  -\frac{n \int_{\Om} h' dV}{|\ga|} <0 .$$
}
\end{remark}
Recently, Ciraolo-Roncoroni \cite{CR} obtained the radial symmetry of the solution to a partially overdetermined problem inside a convex cone in $\SS_{+}^n$, considering the equation $\Delta u + nu = -n$. In the following, we shall consider the same problem in $\SS^n$ without the assumption that $\Omega$ is contained in $\SS_{+}^n$. Instead, we add an assumption that $\Omega$ is a star-shaped domain with respect to the pole $p$. A domain $\Omega \subset \mathbb{S}^n$ is called {\it star-shaped with respect to $p$} if each component of the boundary $\partial \Omega$ can be written as a graph over a geodesic sphere centered at $p$. Now consider the unit sphere $\mathbb{S}^n=I \times \mathbb{S}^{n-1}$ with the warped product metric $g=dr^2 + h(r)^2 g_{\mathbb{S}^{n-1}}$ as before. Note that the interval $I$ is given by $I=[0, \pi)$, which is different from the hemisphere case. For the pole $p$ of the model space and a convex domain $\omega\subset \mathbb{S}^{n-1}$, we can define a convex cone $\mathcal{C}\subset \mathbb{S}^n$ with vertex at $p$ in the same manner:
$$\mathcal{C}=\mathcal{C}_p (\omega):= \{tx : x\in \omega, t \in I\}.$$
Geometrically, $\mathcal{C}=\mathcal{C}_p (\omega)$ is the set of all the unique great semicircles from $p$ to $-p$ passing through $x$ for any $x \in \omega$. Thus, given a convex domain $\omega \subset \mathbb{S}^{n-1}$, it follows that the cone with vertex $-p$ coincides with the cone with vertex at $p$, i.e.,
$$\mathcal{C}_{-p} (\omega) = \mathcal{C}_{p} (\omega).$$
Adopting the $P$-function method used in \cite{CV2019, QX, Weinberger}, we are able to prove the following theorem, which can be seen as a generalization of the results by Ciraolo-Roncoroni \cite{CR} and Pacella-Tralli \cite{PT}.

    \begin{thm}\label{thm2}
        Let $\C \subset \SS^n$ be an open convex cone with vertex at $p$ and $\Om\subset \C$ be an admissible interior domain.  Assume that $\Om$ is a star-shaped domain with respect to $p$ and $-p \notin \partial \Omega$.  Suppose that there exists a solution $u\in C^2(\overline{\Om})$ satisfying
         $$\begin{cases}
            \Delta u +nu = -n & \mbox{in } \Om, \\
            u=0, \quad \frac{\pa u}{\pa \nu} = const = c  & \mbox{on } \ga, \\
            \frac{\pa u}{\pa \nu} = 0 & \mbox{on } \gacone \setminus \{p\},
            \end{cases}$$
where $\nu$ is the outward unit normal to $\ga \cup \gacone \setminus\{p\}.$ Then $\Om$ is part of the geodesic ball $B_R (p_0)$ of radius $R$ centered at $p_0$ in the cone $\mathcal{C}$, i.e.,
            $$\Om = \C \cap B_R (p_0)$$
 and the solution $u$ is given by
            $$u(x) = \frac{1}{\cos R} (\cos r(x) -\cos R)$$
         with $r(x) = \mathrm{dist}(p_0,x)$ and $R=\tan^{-1}(-c)$. Moreover, one of the following two possibilities holds: \\
 (I) $p_0=p$; \\
 (II) $p_0 \in \partial \mathcal{C}$ and $\partial \Omega \cap \partial \mathcal{C}$ is totally geodesic.
    \end{thm}

\begin{proof}
It is well-known that
        \begin{align}
         (\Delta u)^2 \leq n\, \trhess. \label{2.2:2}
        \end{align}
       Note that equality in (\ref{2.2:2}) holds if and only if $\hess u$ is proportional to the metric $g$. Using the Bochner formula,
            \begin{align}\label{boch::star-shaped}
                \Delta |\grad u|^2 & = 2 \<\grad (\Delta u), \grad u\> + 2 \trhess  + 2 \ric (\grad u, \grad u)\nonumber \\
                \nonumber & \geq -2n |\grad u|^2 + \frac{2}{n} (\Delta u)^2 + 2(n-1) |\grad u|^2 \\
                \nonumber& = \frac{2}{n} (-n -nu) \Delta u - 2 |\grad u|^2 \\
                \nonumber& = -2 \Delta u -2u \Delta u - 2 |\grad u|^2\\
                 & =-2\Delta u - \Delta u^2.
            \end{align}
        Define the function $P$ by
        $$P(u):= |\grad u|^2 + 2u + u^2.$$
       Then (\ref{boch::star-shaped}) implies
            $$\Delta P \geq 0.$$
      We also define another function $\p$ by
      $$\p := \<\grad u, \grad \f\> + u\f +\f,$$
      where $h(r)= \sin r$ and $r(x) = \mathrm{dist}(p,x)$. Then
            $$\hess h' = -h' g \ \ \mbox{ and }\ \  \Delta h' = -nh',$$
        where $g$ is the metric of $\SS^n$.
        By the polarized Bochner formula, we get
            \begin{align*}
                \Delta &\<\grad u, \grad \f\> \\
                & = \< \grad (\Delta u), \grad \f \> + \< \grad u, \grad (\Delta \f)\> + 2 \tr ( \hess u \circ \hess \f) + 2 \ric (\grad u, \grad \f) \\
                & = -n \< \grad u, \grad \f\> - n \< \grad u,\grad \f\> - 2\f \Delta u + 2(n-1)\<\grad u, \grad \f \> \\
                & = -2\f (-nu-n) -2 \<\grad u, \grad \f\> \\
                & = 2nu\f +2n\f -2 \<\grad u,\grad \f\>.
            \end{align*}
        Since $$\Delta(u\f) = u \Delta \f + \f \Delta u + 2\<\grad u, \grad \f\> = -2nu\f - n\f +2 \<\grad u, \grad \f\>,$$ we obtain
            $$\Delta ( \<\grad u, \grad \f\> + u\f ) = n\f = - \Delta \f,$$
        which shows
        $$\Delta \p =0.$$
        Note that
        $$\frac{\pa u}{\pa \nu} = 0 \ \ \mbox{ and } \ \ \frac{\pa \f}{\pa \nu} = 0 \ \ \mbox{ on } \ \  \gacone,$$
        which implies that $\grad u$ and $\grad \f$ are tangent vectors on $\gacone$. On $\Gamma_1$,
            $$\grad_{\grad \f} \nu = 0 \ \  \mbox{ and } \ \ \frac{\partial u}{\partial \nu}\equiv const.$$
        Thus
            \begin{equation}\label{2.2: 13}
              0 = \grad_{\grad \f} \(\frac{\pa u}{\pa \nu}\) = \< \grad_{\grad \f} \grad u, \nu \> + \< \grad  u, \grad_{\grad \f} \nu \> = \hess u (\grad \f, \nu )
            \end{equation}
        on $\gacone.$ Moreover, by the convexity of the cone $\mathcal{C}$, we have
        $$\mathrm{II} (\grad u, \grad u) \geq 0 \quad  \mbox{ on }  \gacone,$$
        where $\mathrm{II}(\cdot, \cdot)$ is the  second fundamental form. Thus
            \begin{align*}
              0 = \grad_{\grad u} \(\frac{\pa u}{\pa \nu}\) & = \<\grad_{\grad u}\grad u, \nu\> + \< \grad u, \grad_{\grad u}\nu\> \\
               & = \hess u (\grad u, \nu ) + \mathrm{II}(\grad u, \grad u) \\
               &  \geq \hess u(\grad u, \nu) \quad \mbox{ on }  \Gamma_1.
            \end{align*}
        Hence we obtain
            \begin{align*}
              \frac{\pa P}{\pa \nu} &  = 2 \hess u (\grad u, \nu ) + 2 \frac{ \pa u}{\pa \nu} + 2 u \frac{\pa u}{\pa \nu} \leq 0 \quad \mbox{ on }\gacone.
            \end{align*}
        Suppose that neither $P$ nor $\p$ is constant.  We claim that $P \leq c^2$ in $\Om$. To see this, we note that $P$ satisfies the following.
            $$\begin{cases}
                \Delta P \geq 0 & \mbox{in } \Om, \\
                P \equiv c^2 & \mbox{on } \ga, \\
                \frac{\pa P}{\pa \nu} \leq 0 & \mbox{on } \gacone.
            \end{cases}$$
        Using the divergence theorem, we get
            \begin{align*}
                0 & \leq \int_{\Om} (P -c^2)^{+} \Delta P dV \\
                &  = \int_{\Om} \mathrm{div}((P-c^2)^{+} \grad P) dV - \int_{\Om} \< \grad (P-c^2)^{+},\grad P\> dV \\
                & = \int_{\pa \Om} (P-c^2)^{+} \frac{\pa P}{\pa \nu} d\sigma - \int_{\Om \cap \{P >c^2\}} |\grad P|^2 dV \leq 0,
            \end{align*}
        where $(P-c^2)^{+}=\max \{ P-c^2, 0\}$. Thus we see that  $P \leq c^2$ in $\Om$ and $P$ attains its maximum value on $\Gamma_0$.

        Let $\{e_1, \ldots, e_{n-1}, \nu\}$ be a local orthonormal frame of $\Omega$ at $ x \in \ga$, where $\{e_i\}_{i=1}^{n-1}$ is tangent to $\Gamma_0$ and $\nu$ is  orthogonal to $\Gamma_0$. Since $\ga$ is a level set of $u$, we have
            $$u_i=0\quad  \mbox{ and } \quad u_{ij}=0 \quad \mbox{ on }\Gamma_0$$
         for all $i,j=1,\cdots, n-1$. Moreover, since $\frac{\pa u}{\pa \nu}$ is constant on $\Gamma_0$,
         $$u_{\nu i} = 0\quad \mbox{ on } \Gamma_0$$
         for all $i=1,\cdots, n-1$. Thus
            \begin{eqnarray}\label{2.2:22}
                \hess u (\grad u, \nu) = u_{\nu\nu} \frac{\pa u}{\pa \nu} \quad \mbox{ and } \quad
                \hess u (\grad h',\nu) = u_{\nu\nu} \frac{\pa h'}{\pa \nu} \quad \mbox{ on } \Gamma_0.
            \end{eqnarray}
        Applying the Hopf boundary point lemma on $\Gamma_0$,
            $$0 < \frac{\pa P}{\pa \nu} = 2 \hess u (\grad u, \nu) + 2 \frac{\pa u}{\pa \nu} + 2u \frac{\pa u}{\pa \nu} = 2 \frac{\pa u}{\pa \nu} (u_{\nu\nu} + 1) = 2 c (u_{\nu \nu} + 1).$$
        Since $c$ is constant, we obtain
            \begin{equation}\label{unn}
                u_{\nu\nu} + 1 > 0 \quad \mbox{ or } \quad u_{\nu\nu} + 1 < 0 \quad \mbox{ on }\ga.
            \end{equation}
       Note that $$\<\grad \f, \nu\> = -\sin r \< \grad r, \nu\> = 0 \quad \mbox{ on } \gacone$$
        and
            $$\frac{\pa \p}{\pa \nu} = \hess u (\grad \f, \nu ) + \hess \f (\grad u, \nu ) + u \frac{\pa \f}{\pa \nu} + \f \frac{\pa u}{\pa \nu} + \frac{\pa \f}{\pa \nu}. $$
        Thus it follows from (\ref{2.2: 13}) and the boundary condition that $$\frac{\pa \p}{\pa \nu} =0 \quad \mbox{ on } \gacone.$$
        On the other hand, it follows from (\ref{2.2:22}) that
            \begin{align*}
                \frac{\pa \p}{\pa \nu} & = \hess u (\grad \f, \nu ) + \hess \f (\grad u, \nu ) + u \frac{\pa \f}{\pa \nu} + \f \frac{\pa u}{\pa \nu} + \frac{\pa \f}{\pa \nu} \\
                & = u_{\nu\nu} \frac{\pa \f}{\pa\nu} - \f \frac{\pa u}{\pa \nu} + u \frac{\pa \f}{\pa \nu} + \f\frac{\pa u}{\pa \nu} + \frac{\pa \f}{\pa \nu} \\
                & = (u_{\nu\nu} +1) \frac{\pa \f}{\pa \nu} \quad \mbox{ on } \ga.
            \end{align*}
        Since $\Om$ is a star-shaped domain with respect to $p$, we have
            $$\frac{\pa \f}{\pa \nu} = - \sin r \< \grad r, \nu \> <0 \quad \mbox{ on } \ga.$$
        From (\ref{unn}) we deduce that $\frac{\pa \p}{\pa \nu}<0$ or $\frac{\pa \p}{\pa \nu}>0$ on $\ga$. Applying the divergence theorem, we have
            $$0 = \int_{\Om} \Delta \p dV = \int_{\pa \Om} \frac{\pa \p}{\pa \nu} d \sigma = \int_{\ga} \frac{\pa \p}{\pa \nu} d \sigma <0 \quad ( \mbox{or } >0),$$
        which is a contradiction. Therefore either $P$ or $\p$ is a constant function in $\Om$.

    Suppose $\p$ is a constant function. Then $$\frac{\pa \p}{\pa \nu} =0 \quad \mbox{ and }\quad u_{\nu\nu} +1 =0 \quad \mbox{ on }\ga,$$
         which implies that $$\frac{\pa P}{\pa \nu} =0 \quad \mbox{ on } \ga.$$ Since $P$ has the maximum value on $\ga$,  $P$ is constant in $\Om$ by the Hopf boundary point lemma. Thus we may assume that $P$ is a constant function in $\Omega$. In particular, $\Delta P = 0$ in $\Omega$. We see that equality holds in (\ref{boch::star-shaped}), which implies that $\hess u$ is proportional to the metric $g$. Thus
            \begin{equation}\label{2:hessu:temp}
                \hess u = \frac{\Delta u}{n} g = (-u-1) g.
            \end{equation}
       Since $u=0$ on $\overline{\ga}$, the function $u$ defined on $\overline{\Omega}$ cannot attain simultaneously both its maximum and minimum values on $\overline{\ga}$, which shows that $u$ attains either its maximum or minimum value at some point $p_0 \in \Om\cup \gacone$. Then we have the following three cases:  \\

       \noindent
        {\bf (a)} Suppose $p_0 \in \Om$. Clearly, $\grad u(p_0) = 0$. \\
        {\bf (b)} Suppose $p_0 \in \gacone \setminus \{p\}$. Then $\grad u(p_0)$ also vanishes. To see this, note that since $p_0$ is the maximum or minimum point of $u$ on $\overline{\Omega}$, the restriction of $u$ on $\Gamma_1$ has its maximum or minimum value at $p_0$. This shows that $\nabla^T u (p_0)=0,$ where $\nabla^T$ denotes the induced connection on the tangent bundle $T\Gamma_1$. Thus
         \begin{align*}
         \nabla u(p_0) &= \nabla^T u (p_0) + \nabla^\perp u (p_0) \\
         &= \nabla^T u (p_0) + \frac{\pa u}{\pa \nu}(p_0) \\
         &=0,
         \end{align*}
        where $\nabla^\perp$ denotes the induced connection on the normal bundle $N\Gamma_1$.\\
        {\bf (c)} Suppose $p_0=p \in \gacone$. If $\partial \mathcal{C}$ is smooth everywhere, i.e., totally geodesic, then we have $\grad u(p_0) = 0$ by the same reason as in case {\bf (b)}. If $p$ is the singular point of the cone $\mathcal{C}$, then  $\grad u$ cannot be a nonzero vector at $p$ because $\grad u$ is a continuous vector field and $\frac{\pa u}{\pa \nu} =0$ along $\gacone$.\\

      \noindent   Hence, for any case, at the maximum or minimum point $p_0 \in \Omega \cup \Gamma_1$, we have
         $$\grad u(p_0) =0.$$
         Let $\gamma(s)$ be a unit-speed geodesic passing through $p_0$ satisfying
         $$\gamma (0) = p_0, \grad_{\gamma'(s)}\gamma'(s) = 0, \mbox{ and } |\gamma'(s)|^2 = 1.$$
         Let $f(s):= u(\gamma(s))$. Then
            $$f'(s) = \<\grad u, \gamma'(s)\>$$
        and by (\ref{2:hessu:temp})
            \begin{align*}
                f''(s) & = \< \grad_{\gamma'(s)} \grad u, \gamma'(s)\> + \< \grad u, \grad_{\gamma'(s)}\gamma'(s)\> \\
                & = \hess u (\gamma'(s),\gamma'(s)) \\
                & = -1 -f(s).
            \end{align*}
     From the fact that $\grad u(p_0) = 0$, we obtain an initial value problem:
            \begin{align*}
                f''(s) +f(s)  = -1, \quad     f'(0)  = 0.
            \end{align*}
        A general solution to this ODE is given by
            $$f(s) = c_1 \cos s + c_2 \sin s -1,$$
        where $c_1$ and $c_2$ are constants. Using the initial condition,
            $$f(s) = c_1 \cos s -1,$$
        which shows that the solution $u$ depends only on the geodesic distance because $\gamma$ was arbitrarily chosen to be a geodesic passing through $p_0$. Therefore
            \begin{equation}\label{temp}
                u(x) = c_1 \cos r(x) -1,
            \end{equation}
        where $r(x) = \mathrm{dist} (p_0,x)$. Since $u=0$ on $\ga$ and $\cos r$ is injective in $r \in [0, \pi)$, $\ga$ is part of the geodesic sphere centered at $p_0$ with radius $R=\cos^{-1}\(\frac{1}{c_1}\)$.  Since $\Om$ is connected,
            $$u(x) = \frac{1}{\cos R} \(\cos r - \cos R\)$$
        and $\Om= \C\cap B_R (p_0)$, where $B_R (p_0)$ denotes the geodesic ball centered at $p_0$ with radius $R$.  Since
        $$c= \frac{\pa u}{\pa \nu} = - \frac{\sin R}{\cos R} \<\grad r , \nu\> = - \tan R \quad \mbox{ on } \ga,$$
        we get
        $$R= \tan^{-1}(-c).$$
       Observe that $\nabla u (x)$ is parallel to $\nabla r (x)$ by (\ref{temp}). Moreover, $\nabla u(x)$ lies on the tangent space of $\gacone$ for all $x\in\gacone$ by the boundary condition on $\Gamma_1$. Therefore the point $p_0$ satisfies one of the following two possibilities: \\

    \noindent
       (I) $p_0$ is the vertex $p$. \\
      (II) $p_0\in \partial \mathcal{C}$ and $\partial \Omega \cap \partial \mathcal{C}$ is totally geodesic. \\

    \noindent  For (I), we see that $\Omega = \C\cap B_R (p)$. For (II), $\Omega$ is clearly a half geodesic ball centered at $p_0$ lying over a totally geodesic portion of $\mathcal{C}$.

    \end{proof}

\section{2-dimensional partially overdetermined problems outside a convex cone}

In Section $2$ we studied partially overdetermined PDE problems for a domain inside a convex cone. One may ask whether the similar results as Theorem \ref{thm1} and Theorem \ref{thm2} are still valid for a domain {\it outside} a convex cone. In this section, we give a partial answer to this question.

\begin{df}
Let $\C$ be a convex cone with vertex at $p$ in a space form $M$. A connected bounded open set $\Om\subset M \setminus \overline{\C}$ is an {\it admissible exterior} domain if the boundary $\pa \Om$ satisfies the following.
\begin{itemize}
\item $\pa \Om$ contains the vertex $p$.
\item $\ga:=\pa \Om \setminus \pa \C \neq \emptyset$ is an $(n-1)$-dimensional smooth manifold.
\item $\gacone:= \pa \Om \setminus \overline{\ga} \neq \emptyset$ and $\partial \Gamma_0 = \partial \Gamma_1 \subset \partial \mathcal{C} \setminus \{p\}$  is an $(n-2)$-dimensional smooth manifold.
\item $\mathcal{H}^{n-1} (\Gamma_i) >0$ for $i=0,1$, where $\mathcal{H}^{n-1}$ denotes the $(n-1)$-dimensional Hausdorff measure.
\end{itemize}
\end{df}

\begin{figure}[h!t!b]
\begin{center}
 \resizebox{0.3\textwidth}{!}{%
\begin{tikzpicture}
\draw[name path = C]  plot [smooth, tension=0.6] coordinates { (0,0.89) (-0.5,-1) (2,-2) (4,0) (3,0.89)};
\draw[name path = A,draw=none] (0,0.89) -- (3,0.89);
\tikzfillbetween[of=C and A]{gray, opacity=0.2};
\filldraw[white] (0,0.89)  -- (1.5,-0.8)-- (3,0.89);
\draw (4,2) -- (1.5,-0.8) -- (-1,2);
\coordinate[label={\large $p$}] (p) at (1.5,-1.3);
\coordinate[label={\large $\Om$}] (Om) at (0.3,-0.9);
\coordinate[label={\large $\gacone$}] (gacone) at (2.7,-0.3);
\coordinate[label={\large $\ga$}] (ga) at (2.7,-2.3);
\coordinate[label={\large $\C$}] (C) at (1.6,1.1);
\end{tikzpicture}
}
\caption{An admissible exterior domain $\Om$ outside a convex cone $\mathcal{C}$}\label{sector-like:out}
\end{center}
\end{figure}

\noindent Using the same argument as in the proof of Theorem \ref{thm1}, we consider the following partially overdetermined problem outside a convex cone in a 2-dimensional case.

%
%
    \begin{thm}\label{thm3.1}
        Let $M$ be a $2$-dimensional space form $\RR^2$ or $\SS^2$. Let $\C \subset M$ be an open convex cone with vertex at $p$ and $\Om$ be an admissible exterior domain in $M \setminus \overline{\C}$. If $M = \SS^2$,  assume that either $\Om$ is contained in $\SS_{+}^2$ or $u$ is positive on $\Om$ and assume that $-p \notin \partial \Omega$. Suppose there exists a solution $u\in C^2(\overline{\Om})$ satisfying
            $$\begin{cases}
                \Delta u = -2\f & \mbox{in } \Om, \\
                u = 0, \quad \frac{\pa u}{\pa \nu} = const = c & \mbox{on } \ga, \\
                \frac{\pa u}{\pa \nu} = 0 & \mbox{on } \gacone \setminus \{ p\},
            \end{cases}$$
        where $\nu$ is the outward unit normal to $\ga\cup\gacone \setminus\{p\}$ and the function $h(r)$ is the same as in Section $2$ with $r(x) = \mathrm{dist}(p,x)$. Then
            $$\Om = B_R (p) \setminus \overline{\C},$$
        where $B_R (p)$ denotes the geodesic ball centered at the vertex $p$ with radius $R$ in $M$. In particular, the solution $u$ is given by
            $$u(x) = \begin{cases}
                \frac{R^2 - r^2}{2} & \mbox{in } \RR^2, \\
                \cos r - \cos R & \mbox{in } \SS^2.
                \end{cases}$$
    \end{thm}

    \begin{proof}
        First let us assume that $M = \SS^2$. The proof uses the same argument as in the proof of Theorem \ref{thm1}. However we have a simpler situation in dimension $2$. On $\Gamma_1$, the boundary condition  that $\frac{\partial u}{\partial \nu}=0$ and $\frac{\partial h'}{\partial \nu}=0$ implies
        \begin{itemize}
        \item $\nabla_{\nabla h'} \nu = 0$ and $\nabla_{\nabla u} \nu = 0$,
        \item $ 0 = \grad_{\grad u} \(\frac{\pa u}{\pa \nu}\) = \< \grad_{\grad u} \grad u, \nu \> + \< \grad u, \grad_{\grad u} \nu \> = \hess u (\grad u, \nu)$,
        \item $0 = \grad_{\grad \f} \(\frac{\pa u}{\pa \nu}\) = \< \grad_{\grad \f} \grad u, \nu \> + \< \grad u, \grad_{\grad \f } \nu\> = \hess u(\grad h', \nu)$.
        \end{itemize}
        Using this observation and the argument as in the proof of Theorem \ref{thm1}, we can obtain the conclusion. Furthermore, if $M=\RR^2$, then we can prove Theorem \ref{thm3.1} in the same manner.
  \end{proof}

\noindent  Using the same functions $P$ and $\widetilde{P}$ as in Theorem \ref{thm2}, we obtain a similar Serrin-type symmetry result for domains outside a convex cone as follows.

%
%

    \begin{thm}\label{thm3.2}
        Let $\C \subset \SS^2$ be a convex cone with vertex at the pole $p$ and $\Om$ be an admissible exterior domain in $\SS^2 \setminus \overline{\C}$. Assume that $\Om$ is a star-shaped domain with respect to $p$.  Suppose that there exists a solution $u\in C^2(\overline{\Om})$ satisfying
            $$\begin{cases}
                \Delta u + 2u = -2 & \mbox{in } \Om, \\
                u = 0, \quad \frac{\pa u}{\pa \nu} = const = c  & \mbox{on } \ga, \\
                \frac{\pa u}{\pa \nu} = 0 & \mbox{on } \gacone \setminus \{p\},
            \end{cases}$$
        where $\nu$ is the outward unit normal to $\ga\cup \gacone \setminus\{p\}.$ Then $$\Om = B_R (p_0) \setminus \overline{\C},$$ where $B_R (p_0)$ denotes the geodesic ball centered at $p_0$ with radius $R$ and the solution $u$ is given by
            $$u(x) = \frac{1}{\cos R} (\cos r(x) -\cos R),$$
        where $r(x) = \mathrm{dist}(p_0,x)$. Moreover, one of the following two possibilities holds: \\
        (I) $p_0=p$;\\
        (II) $p_0 \in \partial \mathcal{C}$ and $\partial \Omega \cap \partial \mathcal{C}$ is totally geodesic.
    \end{thm}
    \begin{proof}
        As in the proof of Theorem \ref{thm2}, define two $P$-functions as follows:
            $$P(u)= |\grad u|^2+ 2u + u^2 \quad \mbox{ and } \quad  \p(u) = \<\grad u, \grad \f\> + u\f +\f.$$
        Then
            $$ \Delta P \geq 0 \quad \mbox { and } \quad \Delta \p =0.$$
        Since $\frac{\pa u}{\pa \nu} = 0$ and $\frac{\pa \f}{\pa \nu} = 0$ on $\gacone$, $\grad u$ and $\grad \f$ are tangent vectors of $\gacone$. A direct computation gives
            $$\grad_{\grad \f}\nu = 0 \quad \mbox{ on } \gacone.$$
        Since $\frac{\pa u}{\pa \nu}$ is constant on $\gacone$,
            \begin{equation}\label{ex:gacone:two}
                0 = \grad_{\grad \f} \(\frac{\pa u}{\pa \nu}\) = \< \grad_{\grad \f} \grad u, \nu \> + \< \grad u, \grad_{\grad \f } \nu\> = \hess u(\grad h', \nu) \quad \mbox{ on } \gacone.
            \end{equation}
        We note that $\frac{\pa u}{\pa \nu} =0$ on $\gacone$ implies that $\grad u$ is a radial direction along $\gacone$. This leads to
            \begin{equation*}
                \hess u(\grad u, \nu) = 0 \quad \mbox{ on }\gacone.
            \end{equation*}
        Thus we obtain
            $$\frac{\pa P}{\pa \nu} = 2 \hess u (\grad u, \nu) + 2 \frac{\pa u}{\pa \nu} + 2 u \frac{\pa u}{\pa \nu} = 0 \quad \mbox{ on }\gacone.$$
        Suppose neither $P$ nor $\p$ is constant.
        We claim that $P \leq c^2$ in $\Om$. To see this, we note that $P$ satisfies the following.
            $$\begin{cases}
                \Delta P \geq 0 & \mbox{in } \Om, \\
                P \equiv c^2 & \mbox{on } \ga, \\
                \frac{\pa P}{\pa \nu} = 0 & \mbox{on } \gacone.
            \end{cases}$$
        Using the divergence theorem, we get
            \begin{align*}
                0 & \leq \int_{\Om} (P -c^2)^{+} \Delta P dV \\
                &  = \int_{\Om} \mathrm{div}((P-c^2)^{+} \grad P) dV - \int_{\Om} \< \grad (P-c^2)^{+},\grad P\> dV \\
                & = \int_{\pa \Om} (P-c^2)^{+} \frac{\pa P}{\pa \nu} d\sigma - \int_{\Om \cap \{P >c^2\}} |\grad P|^2 dV \leq  0,
            \end{align*}
        where $(P-c^2)^{+}=\max \{ P-c^2, 0\}$. Thus  $P \leq c^2$ in $\Om$ and $P$ attains its maximum value on $\Gamma_0$.

        Let $\{e_1, \nu\}$ be a local orthonormal frame at $x \in \ga$. Since $\ga$ is  a level set of $u$, we obtain
            $$u_1=0 \quad \mbox{ and } \quad u_{11}=0.$$
        Since $\frac{\pa u}{\pa \nu}$ is constant on   $\ga$, we obtain
        $$u_{\nu 1} = 0 \quad \mbox{ on } \ga.$$
        Then we deduce that on $\ga$,
            \begin{equation}\label{2ex::hess u}
                \hess u (\grad u, \nu) = u_{\nu\nu} \frac{\pa u}{\pa \nu} \quad \mbox{ and }  \quad \hess u (\grad \f, \nu) = u_{\nu\nu} \frac{\pa \f}{\pa \nu}.
            \end{equation}
        By the Hopf boundary lemma, we have
            $$0 < \frac{\pa P}{\pa \nu} = 2 \hess u (\grad u, \nu) + 2 \frac{\pa u}{\pa \nu} + 2u \frac{\pa u}{\pa \nu} = 2 u_{\nu\nu}\frac{\pa u}{\pa \nu} + 2\frac{\pa u}{\pa \nu} = 2 c (u_{\nu \nu} + 1)$$
        on $\ga$.
        Since $c$ is constant, we obtain
            \begin{equation}\label{ex::unn}
                u_{\nu\nu} + 1 > 0 \quad\mbox{ or }\quad u_{\nu\nu} + 1<0\quad \mbox{ on }\ga.
            \end{equation}
        We note that $\<\grad \f, \nu\> = -\sin r \<\grad r, \nu\> =0$ on $\gacone$. Moreover
            $$\frac{\pa \p}{\pa \nu} = \hess u (\grad \f, \nu ) + \hess \f (\grad u, \nu ) + u \frac{\pa \f}{\pa \nu} + \f \frac{\pa u}{\pa \nu} + \frac{\pa \f}{\pa \nu}. $$
        Thus it follows from  (\ref{ex:gacone:two}) that
        $$\frac{\pa \p}{\pa \nu} =0 \quad \mbox{ on }\gacone.$$
        By the boundary conditions on $\ga$ and (\ref{2ex::hess u}), we get
            \begin{align*}
                \frac{\pa \p}{\pa \nu} & = \hess u (\grad \f, \nu ) + \hess \f (\grad u, \nu ) + u \frac{\pa \f}{\pa \nu} + \f \frac{\pa u}{\pa \nu} + \frac{\pa \f}{\pa \nu} \\
                & = u_{\nu\nu} \frac{\pa \f}{\pa\nu} - \f \frac{\pa u}{\pa \nu} + u \frac{\pa \f}{\pa \nu} + \f\frac{\pa u}{\pa \nu} + \frac{\pa \f}{\pa \nu} \\
                & = (u_{\nu\nu} +1) \frac{\pa \f}{\pa \nu} \quad \mbox{ on } \ga.
            \end{align*}
        Since $\Om$ is a star-shaped domain  with respect to $p$,
        $$\frac{\pa \f}{\pa \nu} = - \sin r \< \grad r, \nu \> <0\quad \mbox{ on }\ga.$$
        Then we deduce that
            $$\frac{\pa \p}{\pa \nu}<0\quad\mbox{ or } \quad\frac{\pa \p}{\pa \nu}>0\quad\mbox{ on }\ga$$
        from (\ref{ex::unn}). Using the divergence theorem, we have
            $$0 = \int_{\Om} \Delta \p dV = \int_{\pa \Om} \frac{\pa \p}{\pa \nu} d \sigma = \int_{\ga} \frac{\pa \p}{\pa \nu} d \sigma <0 \quad ( \mbox{or } >0),$$
        which is a contradiction.

        Thus either $P$ or $\p$ is a constant function in $\Om$.
        Suppose $\p$ is a constant function. Then $$\frac{\pa \p}{\pa \nu} =0 \quad \mbox{ and }\quad u_{\nu\nu} +1 =0\quad  \mbox{ on }\ga,$$ which implies that
        $$\frac{\pa P}{\pa \nu} =0\quad\mbox{ on } \ga.$$
        Since $P$ has the maximum value on $\ga$, $P$ is a constant function in $\Om$ by the Hopf boundary point lemma. Therefore we may assume that $P$ is constant. Using the same argument as in the proof of Theorem \ref{thm2}, we can show that the solution $u$ is radially symmetric with respect to some point $p_0$ and it is given by
        $$u(x) = \frac{1}{\cos R} \(\cos r(x) - \cos R\),$$
        where $r(x) = \dist (p_0, x)$. Moreover $\Om$ is the intersection of $\C$ and the geodesic ball $B_R (p_0)$ of radius $R$ centered at $p_0$.
    \end{proof}

\section{An eigenvalue problem with mixed boundary conditions in a cone.}

%
%

Let $M^n (K)$ be one of the space forms $\mathbb{R}^n, \mathbb{S}^n_+,$ and $\mathbb{H}^n$ of constant sectional curvature $K=0, 1,$ and $-1$, respectively. Given the Dirichlet eigenvalue problem for a bounded domain $\Omega \subset M^n (K)$
\begin{eqnarray*}
\begin{cases}
\Delta u +\lambda u = 0 & \mbox{in }  \Om, \\
u =0  & \mbox{on } \partial \Omega,
\end{cases}
\end{eqnarray*}
it is known that the following {\it Rellich identity} holds (see \cite{Rellich} for $K=0$ and \cite{Molzon} for $K=1$ or $-1$):
When $K = 0$,
$$\lambda = -\frac{\int_{\partial \Omega} \frac{\pa f}{\pa \nu}\(\frac{\pa u}{\pa \nu}\)^2  d\sigma }{2 \int_{\Om} u^2 dV}$$
and when $K=1$ or $-1$,
$$\lambda = \frac{-n(n-2)K}{4} -\frac{1}{2K}\frac{\int_{\partial \Omega} \frac{\pa f}{\pa \nu} \(\frac{\pa u}{\pa \nu}\)^2 d \sigma}{\int_{\Omega} fu^2 dV}.$$
Here the function $f(r)$ is defined by
        $$f(r) = \begin{cases}
                       -\frac{r^2}{2} & \mbox{if } K=0, \\
                       \cos r &  \mbox{if } K=1, \\
                       \cosh r & \mbox{if } K=-1,
                     \end{cases}$$
where $r(x)=\dist (p,x )$. Motivated by this, we prove an analogue for an eigenvalue problem with mixed boundary conditions for domains inside a (not necessarily convex) cone in the following.

\begin{thm}\label{thm4.1}
Let $M^n (K)$ be one of space forms $\mathbb{R}^n, \mathbb{S}^n_+,$ and $\mathbb{H}^n$ of constant sectional curvature $K=0, 1,$ and $-1$, respectively. Let $\C \subset M$ be a cone with vertex at $p$ and let $\Omega \subset \C$ be an admissible interior domain.  Suppose there exists a function $u\in C^2(\overline{\Om})$ such that
            $$\begin{cases}
            \Delta u +\lambda u = 0 & \mbox{in } \Omega, \\
            u=0 & \mbox{on } \ga, \\
            \frac{\pa u}{\pa \nu} = 0 &  \mbox{on } \gacone \setminus \{p\}.
            \end{cases}$$
        If $K = 0$, then
            $$\lambda = -\frac{\int_{\ga} \frac{\pa f}{\pa \nu}\(\frac{\pa u}{\pa \nu}\)^2  d\sigma }{2 \int_{\Om} u^2 dV}.$$
        If $K=1$ or $-1$, then $$\lambda = \frac{-n(n-2)K}{4} -\frac{1}{2K}\frac{\int_{\ga} \frac{\pa f}{\pa \nu} \(\frac{\pa u}{\pa \nu}\)^2 d \sigma}{\int_{\Omega} fu^2 dV}.$$
        Here the function $f(r)$ is defined as above.
    \end{thm}

    \begin{proof}
        We first prove the case where $K=1$ or $-1$. The function $f$ satisfies
        $$\Delta f = -nKf  \quad \mbox{ and } \quad \hess f = -Kf g,$$
         where $g$ denotes the metric on $M$. By the polarized Bochner formula, we get
            \begin{align*}
                \Delta &\< \grad u, \grad f\> \\
                 &= \< \grad (\Delta u) , \grad f \> + \< \grad u, \grad (\Delta f)\> +  2 \mathrm{tr} (\hess u
                \circ \hess f) + 2 \ric(\grad u, \grad f) \\
                & = - \lambda \< \grad u, \grad f\> -nK \< \grad u, \grad f\> -2Kf \Delta u + 2(n-1)K \<\grad u, \grad f\> \\
                &  = (-\lambda + nK -2K) \<\grad u, \grad f \> + 2K \lambda f u.
            \end{align*}
        Thus
            \begin{align}\label{1}
                \nonumber\int_{\Om} u \Delta \<\grad u, \grad f \> dV = &(- \lambda + nK -2K ) \int_{\Om} u \< \grad u, \grad f\> dV \\&+ 2K \lambda \int_{\Om} fu^2 dV.
            \end{align}
        On $\ga$, $\grad u$ is parallel to $\nu$ because $\ga$ is a level set of $u$. Using Green's identity,
            \begin{align*}
                \int_{\Om} &u \Delta \<\grad u, \grad f \> dV \\
                 & = \int_{\Om} \<\grad u, \grad f\> \Delta u \, dV + \int_{\pa \Om} u \frac{\pa }{\pa \nu} \<\grad u, \grad f\> d\sigma - \int_{\pa \Omega} \<\grad u, \grad f\> \frac{\pa u}{\pa \nu} \,d \sigma \nonumber \\
                & = - \lambda \int_{\Om} u \<\grad u, \grad f\> dV + \int_{\gacone} u \frac{\pa }{\pa \nu} \< \grad u, \grad f\> d\sigma - \int_{\ga} \frac{\pa f}{\pa \nu} \(\frac{\pa u}{\pa \nu}\)^2 d\sigma.
            \end{align*}
        Moreover we note that $\grad_{\grad f} \nu\ = 0$ on $\gacone$. Since $\frac{\pa u}{\pa \nu}$ is constant on $\gacone$ and $\grad f$ is a tangent vector at $x \in \gacone$, we see that
            $$0 = \grad_{\grad f} \(\frac{\pa u}{\pa \nu}\) = \<\grad_{\grad f} \grad u, \nu \> + \< \grad u, \grad_{\grad f} \nu\> = \hess \,u (\grad f, \nu)\quad \mbox{ on }\gacone,$$
        which implies that
            $$\frac{\pa}{\pa \nu} \< \grad u, \grad f\> = \hess \, u(\grad f,\nu) + \hess \, f(\nu, \grad u) = -Kf \frac{\pa u}{\pa \nu} = 0 \quad \mbox{ on } \gacone.$$
        Thus
            \begin{equation}\label{2}
                \int_{\Om} u \Delta \< \grad u, \grad f\> dV = - \lambda \int_{\Om} u \<\grad u , \grad f\> dV - \int_{\ga} \frac{\pa f}{\pa \nu} \(\frac{\pa u}{\pa \nu}\)^2 d\sigma.
            \end{equation}
        By (\ref{1}) and (\ref{2}), we obtain
            \begin{equation}\label{3}
                (n -2)K \int_{\Om} u \< \grad u, \grad f\> dV + 2K \lambda \int_{\Om} f u^2 dV = - \int_{\ga} \frac{\pa f}{\pa \nu} \(\frac{\pa u}{\pa \nu}\)^2 d \sigma.
            \end{equation}
        Applying the divergence theorem,
            \begin{align*}
                \int_{\Om} u \< \grad u, \grad f\> dV & = \frac{1}{2} \int_{\Om} \<\grad (u^2) , \grad f\> dV \\
                & =\frac{1}{2} \int_{\Om} \mathrm{div} (u^2 \grad f) dV - \frac{1}{2} \int_{\Om} u^2 \Delta f dV \\
                & =\frac{1}{2} \int_{\pa \Om}  u^2  \frac{\pa f}{\pa \nu} \, d\sigma - \frac{1}{2} \int_{\Om} u^2 (-nKf) dV \\
                & = \frac{nK}{2} \int_{\Om}  f u^2 dV.
            \end{align*}
        Plugging the above equality into (\ref{3}), we have
            $$\frac{(n-2)n K^2}{2} \int_{\Om} f  u^2  dV + 2 K\lambda  \int_{\Om} fu^2 dV = - \int_{\ga} \frac{\pa f}{\pa \nu} \(\frac{\pa u}{\pa \nu}\)^2 d\sigma.$$
        Therefore
            $$\lambda = \frac{-n(n-2)K}{4} -\frac{1}{2K}\frac{\int_{\ga} \frac{\pa f}{\pa \nu} \(\frac{\pa u}{\pa \nu}\)^2 d \sigma}{\int_{\Omega} fu^2 dV}.$$

        \noindent Now let us consider the case where $K=0$. In this case,
        $$\Delta f = -n \quad \mbox{ and } \quad \hess f = -Id.$$
        Applying the same argument as in the above gives the conclusion.
    \end{proof}

\vskip 0.3cm
\noindent
{\bf Acknowledgment:}
This work was supported by the National Research Foundation of Korea (NRF-2016R1C1B2009778).

\vskip 1cm
\noindent Jihye Lee\\
Department of Mathematics\\
Sookmyung Women's University\\
Cheongpa-ro 47-gil 100, Yongsan-ku, Seoul, 04310, Korea \\
{\tt E-mail:jihye@sookmyung.ac.kr} \\

\bigskip
\noindent Keomkyo Seo\\
Department of Mathematics and Research Institute of Natural Sciences\\
Sookmyung Women's University\\
Cheongpa-ro 47-gil 100, Yongsan-ku, Seoul, 04310, Korea \\
{\tt E-mail:kseo@sookmyung.ac.kr}\\
URL: http://sites.google.com/site/keomkyo/
\end{document}